\documentclass[leqno,12pt]{article}

\usepackage[nohug]{diagrams}
\diagramstyle[labelstyle=\scriptstyle]

\usepackage{latexsym}
\usepackage{amssymb}
\usepackage{amsthm}
\usepackage{amsmath}
\pdfoutput=1 
\newcommand{\nc}{\newcommand}
\nc{\nt}{\newtheorem} 
\nt{thm}{Theorem}[section]
\newtheorem*{claim}{Claim}
\nt{cor}[thm]{Corollary} 
\nt{prop}[thm]{Proposition}
\nt{lem}[thm]{Lemma} 
\nt{defn}[thm]{Definition}
\nt{exa}[thm]{Example} 
\nt{rem}[thm]{Remark}
\nt{ass}[thm]{Assumption}
\nt{alg}[thm]{Algorithm}
\nt{con}[thm]{Conjecture} 
\nc{\ip}[2]{\mbox{$\langle #1,#2 \rangle$}} 
\nc{\pf}{\noindent{\bf Proof\ \ }}
\nc{\finpf}{\hfill{$\Box$}\linespace}
\nc{\linespace}{\vspace{\baselineskip} \noindent} 
\nc{\R}{{\bf R}}
\nc{\cl}{\mbox{\rm cl}\,} 
\nc{\cls}{ \mbox{{\scriptsize {\rm cl}}}\,} 
\nc{\conv}{\mbox{\rm conv}} 
\nc{\rb}{\mbox{\rm rb}\,}
\nc{\ri}{\mbox{\rm ri}\,}
\nc{\inter}{\mbox{\rm int}\,}
\nc{\bd}{\mbox{\rm bd}\,}
\nc{\epi}{\mbox{\rm epi}\,}
\nc{\gph}{\mbox{\rm gph}\,}
\nc{\rge}{\mbox{\rm rge}\,}
\nc{\rgel}{\mbox{\rm {\scriptsize rge}}\,}
\nc{\dom}{\mbox{\rm dom}\,}
\nc{\detr}{\mbox{\rm det}\,}
\nc{\para}{\mbox{\rm par}\,}
\nc{\crit}{\mbox{\rm crit}\,}

\newenvironment{myequation}{\begin{equation}}{\end{equation}}

\newenvironment{myeqnarray*}{\begin{eqnarray*}}{\end{eqnarray*}}
\nc{\bmye}{\begin{myequation}} \nc{\emye}{\end{myequation}}

\def\tto{\;{\lower 1pt \hbox{$\rightarrow$}}\kern -12pt
           \hbox{\raise 2.8pt \hbox{$\rightarrow$}}\;}

\usepackage{enumerate}
\usepackage{graphicx}

\usepackage{url}
\begin{document}

\author{D. Drusvyatskiy\thanks{%
    Department of Operations Research and Information Engineering,
    Cornell University,
    Ithaca, New York, USA;
    {\tt dd379@cornell.edu}.
    Work of Dmitriy Drusvyatskiy on this paper has been partially supported by the NDSEG grant from the Department of Defense.
    {\bf Corresponding Author} Tel: 718-865-6367, Fax: 607-255-9129
    }%
  \and
	A.D. Ioffe\thanks{
	Department of Mathematics, 
	Technion-Israel Institute of Technology, 
	Haifa, Israel 32000; 
	{\tt ioffe@math.technion.ac.il}.
	Research supported in part by the US-Israel Binational Scientific Foundation Grant 2008261. 
	} 
	\and
  A.S. Lewis\thanks{%
  School of Operations Research and Information Engineering,
  Cornell University,
  Ithaca, New York, USA;
  {\tt http://people.orie.cornell.edu/{\raise.17ex\hbox{$\scriptstyle\sim$}}aslewis/}.
  Research supported in part by National Science Foundation Grant DMS-0806057 and by the US-Israel Binational Scientific Foundation Grant 2008261.
}}

\title{\Large The dimension of semialgebraic subdifferential graphs}

\maketitle
%\begin{definition}
%Let $f:\mathbb{R}^n\rightarrow \mathbb{R}$. Then we define $$\textrm{graph}(f)=\{(x,y)\in\mathbb{R}^{n}\times\mathbb{R}:y=f(x)\textrm{ and } y<\infty\}$$
%\end{definition}

\begin{abstract}
Examples exist of extended-real-valued closed functions on  ${\bf R}^n$  whose subdifferentials (in the standard, limiting sense) have large graphs.  By contrast, if such a function is semi-algebraic, then its subdifferential graph must have everywhere constant local dimension $n$.    This result is related to a celebrated theorem of Minty, and surprisingly may fail for the Clarke subdifferential.
%We prove that the graph of the limiting subdifferential of any extended-real-valued, lower-semicontinuous, semi-algebraic function defined on $\R^n$, has constant local dimension $n$ at each of its points. In contrast, this may fail for the Clarke subdifferential. We discuss connections between our results and a celebrated theorem of Minty.  
\end{abstract}


\section{Introduction}
A principle goal of variational analysis and nonsmooth optimization (and of critical point theory) is to study generalized critical points of extended-real-valued functions on $\R^n$. These are the points where a generalized subdifferential, such as the Frechet, limiting, or Clarke subdifferential, of $f$ contains the zero vector. Generalized critical points of smooth functions are, in particular, critical points in the classical sense, while critical points of convex functions are simply their minimizers. More generally, one could consider the perturbed function $x\mapsto f(x)-\langle v,x\rangle$, for some fixed vector $v\in\R^n$. Then a point $x$ is critical precisely when the pair $(x,v)$ lies in the graph of the subdifferential. Hence, it is natural to try to understand geometric properties of subdifferential graphs. 

In particular, an interesting question in this area is to understand the ``size'' of the subdifferential graph. For instance, for a smooth function defined on $\R^n$, the graph of the subdifferential is an $n$-dimensional surface. Minty \cite{minty} famously showed that the subdifferential graph of a lower semicontinuous, convex function defined on $\R^n$ is Lipschitz homeomorphic to $\R^n$. In fact, he provided explicit Lipschitz homeomorphisms that are very simple in nature. More generally in \cite{prox_reg}, Poliquin and Rockafellar used Minty's theorem to show that an analogous result holds for ``prox-regular functions'', unifying the smooth and the convex cases. Hence, we would expect that for a nonpathological function, the subdifferential graph should have the same dimension, in some sense, as the space that the function is defined on. A limiting feature of Poliquin's and Rockafellar's approach is that their arguments rely on convexity, or rather the related notion of maximal monotonicity. Hence their techniques do not seem to extend to a larger class of functions.

From a practical point of view, the size of the subdifferential graph may have important algorithmic applications. For instance, Robinson \cite{rob} shows computational promise for functions defined on $\R^n$ whose subdifferential graphs are locally homeomorphic to an open subset of $\R^n$. In particular, due to Minty's result, Robinson's techniques are applicable for lower semicontinuous, convex functions. When can we then be sure that the dimension of the subdifferential graph is the same as the dimension of the domain space?     

It is well-known that for general functions, even ones that are Lipschitz continuous, the subdifferential graph can be very large. For instance, there is a  1-Lipschitz function $f\colon\R\to\R$, such that the Clarke subdifferential $\partial_c f$ is the unit interval $[-1,1]$ at every point. Furthermore, this behavior is typical \cite{gen_lip} and such pathologies are not particular to the Clarke case \cite{Benoist, badlim}. 

These pathological functions, however, do not normally appear in practice. As a result, the authors of \cite{dim} were led to consider {\em semi-algebraic} functions, those functions whose graphs are defined by finitely many polynomial equalities and inequalities. They showed that for a proper, semi-algebraic function on $\R^n$, any reasonable subdifferential has a graph that is, in a precise mathematical sense, exactly $n$-dimensional. The authors derived a variety of applications for generic semi-algebraic optimization problems.
 
The dimension of a semi-algebraic set, as discussed in \cite{dim}, is a global property governed by the maximal size of any part of this set. In particular, the result above does not rule out that some parts of the subdifferential graph may be small. In fact, in the Clarke case this can happen! It is the aim of our current work to elaborate on this phenomenon and to show that it does not occur in the case of the limiting subdifferential. Specifically, we will show that for a lower semicontinuous, semi-algebraic function $f$ on $\R^n$, the graph of the limiting subdifferential has local dimension $n$, uniformly over the whole set. Surprisingly, as we noted, this type of a result does not hold for the Clarke subdifferential. That is, even for the simplest of examples, the graph of the Clarke subdifferential may be small in some places, despite being a larger set than the limiting subdifferential graph.

To be concrete, we state our results for semi-algebraic functions.  Analogous results, with essentially identical proofs, hold for functions definable in an ``o-minimal structure'' and, more generally, for ``tame'' functions. In particular, our results hold for globally subanalytic functions, discussed in \cite{Shiota}. For a quick introduction to these concepts in an optimization context, see \cite{tame_opt}.

\section{Preliminaries}
\subsection{Variational Analysis}
In this section, we summarize some of the fundamental tools used in variational analysis and nonsmooth optimization.
We refer the reader to the monographs Borwein-Zhu \cite{Borwein-Zhu},
Mordukhovich \cite{Mord_1,Mord_2}, Clarke-Ledyaev-Stern-Wolenski \cite{CLSW}, and Rockafellar-Wets \cite{VA}, for more details.  Unless otherwise stated, we follow the terminology and notation of \cite{VA}.

The functions that we will be considering will be allowed to take values in the extended real line $\overline{\R}:=\R\cup\{-\infty\}\cup\{+\infty\}$. We say that an extended-real-valued function is proper if it is never $-\infty$ and is not always $+\infty$.  

For a function $f\colon\R^n\rightarrow\overline{\R}$, we define the {\em domain} of $f$ to be $$\mbox{\rm dom}\, f:=\{x\in\R^n: f(x)<+\infty\},$$ and we define the {\em epigraph} of $f$ to be the set $$\mbox{\rm epi}\, f:= \{(x,r)\in\R^n\times\R: r\geq f(x)\}.$$

A {\em set-valued mapping} $F$ from $\R^n$ to $\R^m$, denoted by $F\colon\R^n\rightrightarrows\R^m$, is a mapping from $\R^n$ to the power set of $\R^m$. Thus for each  point $x\in\R^n$, $F(x)$ is a subset of $\R^m$. For a set-valued mapping $F\colon\R^n\rightrightarrows\R^m$, the {\em domain}, {\em graph}, and {\em range} of $F$ are defined to be $$\mbox{\rm dom}\, F:=\{x\in\R^n:F(x)\neq\emptyset\},$$ $$\mbox{\rm gph}\, F:=\{(x,y)\in\R^n\times\R^m:y\in F(x)\},$$ $$\rge F= \bigcup_{x\in\R^n} F(x),$$ respectively. Observe that $\dom F$ and $\rge F$ are images of $\gph F$ under the projections $(x,y)\mapsto x$ and $(x,y)\mapsto y$, respectively.

Throughout this work, we will only use Euclidean norms. Hence for a point $x\in\R^n$, the symbol $|x|$ will denote the standard Euclidean norm of $x$. Given a point $\bar{x}\in\R^n$, we let $o(|x-\bar{x}|)$ be shorthand for a function that satisfies $\frac{o(|x-\bar{x}|)}{|x-\bar{x}|}\rightarrow 0$ whenever $x\to \bar{x}$ with $x\neq\bar{x}$. We now turn to {\em subdifferentials}, which are fundamental objects in variational analysis.

\begin{defn}
{\rm Consider a function $f\colon\R^n\to\overline{\R}$ and a point $\bar{x}$ with $f(\bar{x})$ finite. 
\begin{enumerate}
\item
The {\em Frechet subdifferential} of $f$ at $\bar{x}$, denoted 
$\hat {\partial} f(\bar x)$, consists of all vectors $v \in \R^n$ such that $$f(x)\geq f(\bar{x})+\langle v,x-\bar{x} \rangle + o(|x-\bar{x}|).$$ 
\item We define the {\em Frechet subjet} of $f$ to be the set $$[\hat{\partial}f]=\{(x,y,v)\in\R^n\times\R\times\R^n:y=f(x), v\in\hat{\partial} f(x)\}.$$
\end{enumerate}}
\end{defn}

The Frechet subjet does not have desirable closure properties. Consequently, the following definition is introduced.

\begin{defn}
{\rm Consider a function $f\colon\R^n\to\overline{\R}$ and a point $\bar{x}$ with $f(\bar{x})$ finite. 
\begin{enumerate}
\item The {\em limiting subdifferential} of $f$ at $\bar{x}$, denoted 
$\partial f(\bar x)$, consists of all vectors $v \in \R^n$ such that there is a sequence $(x_i,f(x_i),v_i)\in[\hat{\partial}f]$ with $(x_i,f(x_i),v_i)\to(\bar{x},f(\bar{x}),v)$.
\item We define the {\em limiting subjet} of $f$ to be the set $$[\partial f]=\{(x,y,v)\in\R^n\times\R\times\R^n:y=f(x), v\in\partial f(x)\}.$$
\end{enumerate}}
\end{defn}

For $x$ such that $f(x)$ is not finite, we follow the convention that $\hat{\partial}f(x)=\partial f(x)=\emptyset$. The following is a standard result in subdifferential calculus.
\begin{prop}\cite[Exercise 10.10]{VA}\label{prop:lip}
Consider a function $f_1\colon\R^n\to\overline{\R}$ that is locally Lipschitz around a point $\bar{x}\in\R^n$ and a function $f_2\colon\R^n\to\overline{\R}$ that is lower semi-continuous and proper with $f_2(\bar{x})$ finite. Then the inclusion$$\partial (f_1+ f_2)(\bar{x})\subset \partial f_1(\bar{x})+\partial f_2(\bar{x}),$$ holds.
\end{prop}

We will have occasion to talk about restrictions of subjets. Given a function $f\colon\R^n\to\overline{\R}$ and a set $M\subset\R^n$, we define the restriction of $[\partial f]$ to $M$ to be the set $[\partial f]\big|_M:=[\partial f]\cap (M\times\R\times\R^n)$. Analogous notation will be used for restrictions of the Frechet subjet $[\hat{\partial}f]$. Observe that in general, the set $[\partial f]\big|_M$
is not a subjet of any function.
%sets $[\partial f]\big|_M$ and $[f\big|_M]$ are in general very different. The former is simply a portion of the subjet $[\partial f]$, while the latter is itself a limiting subjet, but of the restricted function $f\big|_M$. 
More generally, for a set $F\subset\R^n\times\R\times\R^n$ and a set $M\subset\R^n$, we let $F\big|_M:=F\cap (M\times\R\times\R^n)$.

An open ball of radius $r$ around a point $x\in\R^n$ will be denoted by $B_r(x)$, while the closed unit ball of radius $r$ around a point $x\in\R^n$ will be denoted by $\bar{B}_r(x)$. The open and the closed unit balls will be denoted by $\bf{B}$ and $\bf{\overline{B}}$, respectively. Consider a set $M\subset\R^n$. We denote the topological closure, interior, and boundary of $M$ by $\cl M$, $\inter M$, and $\bd M$, respectively. We define the {\em indicator function} of $M$, $\delta_M\colon\R^n\to\overline{\R}$, to be $0$ on $M$ and $+\infty$ elsewhere.  
Indicator functions allow us to translate analytic information about functions to geometric information about sets. In this spirit, we now define {\em normal cones}, which are the geometric analogues of subdifferentials.  

\begin{defn}
{\rm Consider a set $M\subset\R^n$ and a point $x\in\R^n$. The {\em Frechet} and the {\em limiting normal cones} are defined to be $\hat{N}_M(x):=\hat{\partial}\delta_M (x)$ and $N_M(x):=\partial\delta_M (x)$, respectively.
}
\end{defn}

Given any set $Q\subset\R^n$ and a mapping $F\colon Q\to \widetilde{Q}$, where $\widetilde{Q}\subset\R^m$, we say that $F$ is $\bf{C}^1$-{\em smooth} if for each point $\bar{x}\in Q$, there is a neighborhood $U$ of $\bar{x}$ and a $\bf{C}^1$ mapping $\hat{F}\colon \R^n\to\R^m$ that agrees with $F$ on $Q\cap U$. Henceforth, the word {\em smooth} will always mean $\bf{C}^1$-smooth. Since we will not need higher order of smoothness in our work, no ambiguity should arise.
If a smooth function $F$ is bijective and its inverse is also smooth, then we say that $F$ is a {\em diffeomorphism}. More generally, we have the following definition.
\begin{defn}
{\rm
Consider sets $Q\subset\R^n$, $\widetilde{Q}\subset\R^m$, and a mapping $F\colon Q\to \widetilde{Q}$. We say that $F$ is a {\em local diffeomorphism} around a point $\bar{x}\in Q$ if there exists a neighborhood $U$ of $\bar{x}$ such that the restriction 
\begin{equation}
F\big|_{Q\cap U}\colon Q\cap U\to F(Q\cap U),\label{eq:defn}
\end{equation}
is a diffeomorphism. Now consider another set $K\subset\R^m$. We say that $F$ is a local diffeomorphism around $\bar{x}$ onto $K$ if there exists a neighborhood $U$ of $\bar{x}$ such that the mapping in (\ref{eq:defn}) is a diffeomorphism and $K=F(Q\cap U)$.
}
\end{defn}

We now recall the notion of a {\em manifold}.
%\begin{defn}
%For a function $f\colon\R^n\to\overline{R}$, we define the 
%\end{defn}

\begin{defn}[{\cite[Proposition 8.12]{Lee}}]
{\rm Consider a set $M\subset\R^n$. We say that $M$ is a {\em manifold} of dimension $r$ if for each point $\bar{x}\in M$, there is an open neighborhood $U$ around $\bar{x}$ such that $M\cap U=F^{-1}(0)$, where $F\colon U\to\R^{n-r}$ is a $\bf{C}^1$ smooth map with $\nabla F(\bar{x})$ of full rank. In this case, we call $F$ a {\em local defining function} for $M$ around $\bar{x}$.  
}
\end{defn}
Strictly speaking, what we call a manifold is usually referred to as a $\bf{C}^1$-{\em submanifold} of $\R^n$. For a manifold $M\subset\R^n$ and a point $x\in M$, the Frechet normal cone, $\hat{N}_M(x)$, and the limiting normal cone, $N_M(x)$, coincide and are equal to the normal space, in the sense of differential geometry. For more details, see for example \cite[Example 6.8]{VA}. For a smooth map $F\colon M\to N$, where $M$ and $N$ are manifolds, we say that $F$ has {\em constant rank} if its derivative has constant rank throughout $M$.

%We record the following change of coordinates formula for reference.
%\begin{thm}[{\cite[Exercise 6.7]{VA}}]\label{thm:change_coordinates}
%Let $U\subset\R^n$ and $\widetilde{U}\subset\R^m$ be open sets and let $F\colon U\rightarrow\widetilde{U}$ be a smooth map. Let $C=F^{-1}(D)$ for a set $D\subset\widetilde{U}$ and suppose that $\nabla F(\bar{x})$ has full rank $m$ at a point $\bar{x}\in C$. Then $$N_C(\bar{x})=\nabla F(\bar{x})^{*}N_D(F(\bar{x})),$$ $$\hat{N}_C(\bar{x})=\nabla F(\bar{x})^{*}\hat{N}_D(F(\bar{x})).$$ 
%\end{thm}

For a set $M\subset\R^n$ and a point $x\in\R^n$, the distance of $x$ from $M$ is   $$d_M(x)=\inf_{y\in M}|x-y|,$$ and the projection of $x$ onto $M$ is $$P_M(x)=\{y\in M:|x-y|=d_M(x)\}.$$ 
%\begin{defn}
%Consider a mapping $F\colon\R^n\rightrightarrows\R^m$ and a point $\bar{x}\in \mbox{\rm dom}\, F$. The coderivative is the mapping $D^{*}F(\bar{x}|\bar{u})\colon\R^m\rightrightarrows\R^n$ defined by $$ D^{*}F(\bar{x}|\bar{u})(y)=\{v\in\R^n|(v,-y)\in N_{\textrm{graph } F}(\bar{x},\bar{u})\}.$$
%\end{defn}

%\noindent We define $D^*_cF(\bar{x}|\bar{u})$ in an analogous way with the limiting normal cone replaced by the Clarke normal cone. 

%\begin{thm}{\cite[Example 8.53]{VA}}\label{thm:sub_dist}
%For a closed set $M\subset\R^n$, one has at any point $\bar{x}\in S$ that $\hat{\partial}d_M(\bar{x})=\hat{N}_M(\bar{x})\cap \B$.
%\end{thm}
Finally, we will need the following result.
\begin{thm}{\cite[Example 10.32]{VA}}\label{thm:mar}
For a closed set $M\subset\R^n$, the inclusion $$\partial[d_M^2](x)\subset 2[x-P_M(x)],$$ holds for all $x\in\R^n$.
\end{thm}

\subsection{Semi-algebraic geometry}
A {\em semi-algebraic} set $S\subset\R^n$ is a finite union of sets of the form $$\{x\in \R^n: P_1(x)=0,\ldots,P_k(x)=0, Q_1(x)<0,\ldots, Q_l(x)<0\},$$ where $P_1,\ldots,P_k$ and $Q_1,\ldots,Q_l$ are polynomials in $n$ variables. In other words, $S$ is a union of finitely many sets, each defined by finitely many polynomial equalities and inequalities. A map $F\colon\R^n\rightrightarrows\R^m$ is said to be {\em semi-algebraic} if $\mbox{\rm gph}\, F\subset\R^{n+m}$ is a semi-algebraic set. Semi-algebraic sets enjoy many nice structural properties. We discuss some of these properties in this section. For more details, see the monographs of Basu-Pollack-Roy \cite{ARAG}, Lou van den Dries \cite{LVDB}, and Shiota \cite{Shiota}. For a quick survey, see the article of van den Dries-Miller \cite{DM} and the surveys of Coste \cite{Coste-semi, Coste-min}. Unless otherwise stated, we follow the notation of \cite{DM} and \cite{Coste-semi}. 

A fundamental fact about semi-algebraic sets is provided by the Tarski-Seidenberg Theorem \cite[Theorem 2.3]{Coste-semi}. Roughly speaking, it states that a linear projection of semi-algebraic set remains semi-algebraic. From this result, it follows that a great many constructions preserve semi-algebraicity. In particular, for a semi-algebraic function $f\colon\R^n\to\overline{\R}$, it is easy to see that the set-valued mappings $\hat{\partial} f$, $\partial f$, along with the subjets $[\hat{\partial}f]$, $[\partial f]$, are semi-algebraic. See for example \cite[Proposition 3.1]{tame_opt}. 

\begin{defn}
{\rm Given finite collections $\{B_i\}$ and $\{C_j\}$ of subsets of $\R^n$, we say that $\{B_i\}$ is {\em compatible} with $\{C_j\}$ if for all $B_i$ and $C_j$, either $B_i\cap C_j=\emptyset$ or $B_i\subset C_j$.}
\end{defn}

\begin{defn}\label{defn:whit}
{\rm Consider a semi-algebraic set $Q$ in $\R^n$. A {\em stratification} of $Q$ is a finite partition of $Q$ into disjoint, connected, semi-algebraic manifolds $M_i$ (called strata) with the property that for each index $i$, the intersection of the closure of $M_i$ with $Q$ is the union of some $M_j$'s.} 
\end{defn}

The most striking and useful fact about semi-algebraic sets is that stratifications of semi-algebraic sets always exist. In fact, a more general result holds, which is the content of the following theorem. 
\begin{thm}[{\cite[Theorem 4.8]{DM}}]\label{thm:strat}
Consider a semi-algebraic set $S$ in $\R^n$ and a semi-algebraic map $f\colon S\rightarrow\R^m$. 
Then there exists a stratification $\mathcal{A}$ of $S$ and a stratification $\mathcal{B}$ of $\R^m$ such that for every stratum $M\in\mathcal{A}$, we have that the restriction $f|_M$ is smooth, $f(M)\in\mathcal{B}$, and $f$ has constant rank on $M$.
Furthermore, if $\mathcal{A}'$ is some other stratification of $S$, then we can ensure that $\mathcal{A}$ is compatible with $\mathcal{A}'$. 
\end{thm}

\begin{defn}
{\rm Let $A\subset\R^n$ be a nonempty semi-algebraic set. Then we define the {\em dimension} of $A$, $\dim A$, to be the maximal dimension of a stratum in any stratification of $A$. We adopt the convention that $\dim \emptyset=-\infty$.}  
\end{defn}

It can be easily shown that the dimension does not depend on the particular stratification. Dimension is a very well behaved quantity, which is the content of the following proposition. See \cite[Chapter 4]{LVDB} for more details.
\begin{thm}Let $A$ and $B$ be nonempty semi-algebraic sets in $\R^n$. Then the following hold.
\begin{enumerate}
\item If $A\subset B$, then $\dim A\leq \dim B$.
\item $\dim A=\dim \mbox{\rm cl}\,{A}$.
\item $\dim (\mbox{\rm cl}\,{A}\setminus A)< \dim A$.
\item If $f\colon A\rightarrow\R^n$ is a semi-algebraic mapping, then $\dim f(A)\leq \dim A$. If $f$ is one-to-one, then $\dim f(A)=\dim A$. In particular, semi-algebraic homeomorphisms preserve dimension.
\item $\dim A\cup B= \max\{\dim A, \dim B\}$.
\item $\dim A\times B=\dim A+\dim B$.
\end{enumerate}
\end{thm}

Observe that the dimension of a semi-algebraic set only depends on the maximal dimensional manifold in a stratification. Hence, dimension is a somewhat crude measure of the size of the semi-algebraic set. In particular, it does not provide much insight into what the set looks like locally around each of its point. Hence, this motivates a localized notion of dimension.  

\begin{defn}
{\rm Consider a semi-algebraic set $Q\subset \R^n$ and a point $\bar{x}\in Q$. We let the {\em local dimension} of $Q$ at $\bar{x}$ be $$\dim_Q(\bar{x}):=\inf_{r>0}\dim (Q\cap B_r(\bar{x})).$$ In fact, it is not hard to see that there exists a real number $\bar{r}>0$ such that for every real number $0<r<\bar{r}$, we have $\dim_Q(\bar{x})=\dim (Q\cap B_{r}(\bar{x}))$.
}
\end{defn}

The following is now an easy observation.
\begin{prop}\label{prop:dim}\cite[Exercise 3.19]{Coste-semi}
For any semi-algebraic set $Q\subset\R^n$, we have the identity $$\dim Q=\max_{x\in Q}\dim_Q(x).$$
\end{prop}

\begin{defn}~\label{def:triv}
{\rm Let $A\subset\R^m$ be a semi-algebraic set. A continuous semi-algebraic mapping $p\colon A\rightarrow\R^n$ is {\em semi-algebraically trivial} over a semi-algebraic set $C\subset\R^n$ if there is a semi-algebraic set $F$ and a semi-algebraic homeomorphism $h\colon p^{-1}(C)\rightarrow C\times F$ such that $p|_{p^{-1}(C)}={{\rm proj}_C}\circ h$, or in other words the following diagram commutes:}
\begin{diagram}[height=1.7em]
p^{-1}(C) &\rTo^h   &C\times F\\
          &\rdTo_p  &\dTo_{\mbox{\scriptsize {\rm proj$_C$}}} \\
          &         &C
\end{diagram}
{\rm We call $h$ a {\em semi-algebraic trivialization} of $p$ over $C$.}
\end{defn} 
Henceforth, we use the symbol $\cong$ to indicate that two semi-algebraic sets are semi-algebraically homeomorphic.
\begin{rem} \label{rmk:Hardt}
{\rm If $p$ is trivial over some semi-algebraic set $C$, then we can decompose $p|_{p^{-1}(C)}$ into a homeomorphism followed by a simple projection. Also, since the homeomorphism $h$ in the definition is surjective and $p|_{p^{-1}(C)}={\rm proj}_C\circ h$, it easily follows that for any point $c\in C$, we have $p^{-1}(c)\cong F$ and $p^{-1}(C)\cong C\times p^{-1}(c)$.}
\end{rem}

\begin{defn}
{\rm In the notation of Definition~\ref{def:triv}, a trivialization $h$ is {\em compatible} with a semi-algebraic set $B\subset A$ if there is a semi-algebraic set $H\subset F$ such that $h(B\cap p^{-1}(C))= C\times H$.}
\end{defn}

If $h$ is a trivialization over $C$ then, certainly, for any set $B\subset A$ we know $h$ restricts to a homeomorphism from $B\cap p^{-1}(C)$ to $h(B\cap p^{-1}(C))$. The content of the definition above is that if $p$ is compatible with $B$, then $h$ restricts to a homeomorphism between $B\cap p^{-1}(C)$ and the product $C\times H$ for some semi-algebraic set $H\subset F$.
 
The following is a remarkably useful theorem \cite[Theorem 4.1]{Coste-semi}.
\begin{thm}[Hardt triviality]\label{theorem:Hardt}
Let $A\subset\R^n$ be a semi-algebraic set and $p\colon A\rightarrow\R^m$, a continuous semi-algebraic mapping. Then, there is a finite partition of the image $p(A)$ into semi-algebraic sets $C_1,\ldots, C_k$ such that $p$ is semi-algebraically trivial over each $C_i$. Moreover, if $Q$ is a semi-algebraic subset of $A$, we can require each trivialization $h_i\colon p^{-1}(C_i)\rightarrow C_i\times F_i$ to be compatible with $Q$.
\end{thm}

For an application of Hardt triviality to semi-algebraic set-valued analysis, see \cite[Section 2.2]{dim}. 
The following proposition is a simple consequence of Hardt triviality.
\begin{prop}\label{prop:fiber}
Consider semi-algebraic sets $M$ and $Q$ satisfying $M\subset Q\subset\R^n$. Assume that there exists a continuous mapping $p\colon Q\to\R^m$, for some positive integer $m$, such that for each point $x$ in the image $p(Q)$ we have $\dim p^{-1}(x)=\dim (p^{-1}(x)\cap M)$. Then $M$ and $Q$ have the same dimension.
\end{prop}
{\pf Applying Theorem~\ref{theorem:Hardt} to the map $p$, we partition the image $p(Q)$ into finitely many disjoint sets $C_1,\ldots,C_k$ such that for each index $i$, we have the relations
$$p^{-1}(C_i)\cong C_i\times p^{-1}(c),$$
$$p^{-1}(C_i)\cap M \cong C_i\times (p^{-1}(c)\cap M),$$
where $c$ is any point in $C_i$. Since by assumption, the equation $\dim p^{-1}(x)=\dim (p^{-1}(x)\cap M)$ holds for all points $x$ in the image $p(Q)$, we deduce $$\dim p^{-1}(C_i)=\dim (p^{-1}(C_i)\cap M),$$ for each index $i$. Thus 
\begin{align*}
\dim Q=\dim \bigcup_i p^{-1}(C_i)&= \max_i \dim p^{-1}(C_i)=\max_i \dim (p^{-1}(C_i)\cap M)\\
&=\dim \bigcup_i (p^{-1}(C_i)\cap M)=\dim M,
\end{align*}
as we needed to show.
}\qed

We will have occasion to use the following simple proposition \cite[Theorem 3.18]{Coste-min}.
\begin{prop}\label{prop:const_gen}
Consider a semi-algebraic, set-valued mapping $F\colon\R^n\rightrightarrows\R^m$. Suppose there exists an integer $k$ such that the set $F(x)$ is $k$-dimensional for each point $x\in \dom F$. Then the equality, $$\dim\gph F= \dim \dom F +k,$$ holds.   
\end{prop}

%\begin{thm}Let $A$ and $B$ be nonempty semi-algebraic sets in $\R^n$. Then the following hold.
%\begin{enumerate}
%\item If $A\subset B$, then $\dim A\leq \dim B$.
%\item $\dim A=\dim \mbox{\rm cl}\,{A}$.
%\item $\dim (\mbox{\rm cl}\,{A}\setminus A)< \dim A$.
%\item If $f\colon A\rightarrow\R^n$ is a semi-algebraic function, then we have $\dim f(A)\leq \dim A$. If $f$ is one-to-one, then we have $\dim f(A)=\dim A$. In particular, semi-algebraic homeomorphisms preserve dimension.
%\item $\dim A\cup B= \max\{\dim A, \dim B\}$.
%\item $\dim A\times B=\dim A+\dim B$.
%\end{enumerate}
%\end{thm}

\section{Main results}
In our current work, we build on the following theorem. This result and its consequences for generic semi-algebraic optimization problems are discussed extensively in \cite{dim}. 
\begin{thm}\cite[Theorem 3.6]{dim}\label{thm:grd}
Let $f\colon\R^n\rightarrow\overline{\R}$ be a proper semi-algebraic function. Then the graphs of the Frechet and the limiting subdifferentials have dimension exactly $n$.
\end{thm}
In fact, Theorem~\ref{thm:grd} also holds for the proximal and Clarke subdifferentials. For more details see~\cite{dim}.

To motivate our current work, consider a manifold $Q\subset\R^n$. The set, $$\gph N_Q=\{(x,y)\in\R^n\times\R^n:y\in N_Q(x)\},$$ is the normal bundle of $Q$, and as such, $\gph N_Q$ is itself a manifold of dimension $n$ \cite[Proposition 10.18]{Lee}. In particular, $\gph N_Q$ is $n$-dimensional, {\em locally} around each of its points. This suggests that perhaps Theorem \ref{thm:grd} may be strengthened to pertain to the local dimension of the graph of the subdifferential. Indeed, this is the case. In fact, we will prove something stronger.

Let $f\colon\R^n\rightarrow\overline{\R}$ be a lower semicontinuous, proper, semi-algebraic function. Observe that the sets $\gph \partial f$ and $[\partial f]$ are in semi-algebraic bijective correspondence, via the map $(x,v)\longmapsto (x,f(x),v)$, and hence these two sets have the same dimension. Thus by Theorem~\ref{thm:grd}, the dimension of the subjet $[\partial f]$ is exactly $n$. Combining this observation with Proposition~\ref{prop:dim}, we deduce that the local dimension of $[\partial f]$ at each of its points is at most $n$. In this work, we prove that, remarkably, the local dimension of $[\partial f]$ at each of its points is exactly $n$ (Theorem~\ref{thm:local_dim}). From this result, it easily follows that the local dimension of $\gph\partial f$ at each of its points is exactly $n$ as well. Analogous result holds for the Frechet subjet $[\hat{\partial}f]$. 

The proof of Theorem~\ref{thm:local_dim} relies on a very general accessibility result, which we establish in Lemma~\ref{lem:main}. This result, in fact, holds in the absence of semi-algebraicity. In Remark~\ref{rem:fail}, we provide a simple example illustrating that the assumption of lower-semicontinuity is necessary for our conclusions to hold. Then in Subsection~\ref{sec:Clarke}, we recall the definition of the Clarke subdifferential mapping and show that its graph may have small local dimension at some of its points. Thus, the analogue of Theorem~\ref{thm:local_dim} fails for the Clarke subdifferential. This further illustrates the subtlety involved when analyzing local dimension.    

\subsection{Geometry of the Frechet and limiting subdifferential mappings}

\begin{lem}[{\rm Accessibility}]\label{lem:main}
Let $f\colon\R^n\to\overline\R$ be a lower semicontinuous function and $M\subset\R^n$ a closed set on which $f$ is finite. Fix a point $\bar{x}\in M$ and consider a triple $(\bar{x},f(\bar{x}),\bar{v})\in[\hat{\partial} f]\big|_{M}$. Suppose that there exists a sequence of real numbers $m_i\to\infty$ such that $$\bar{v}\in \bd \bigcup_{x\in M} \partial (f(\cdot) + \frac{1}{2}m_i|\cdot-\bar{x}|^2)(x),$$ for each $i$. Then the inclusion $(\bar{x},f(\bar{x}),\bar{v})\in\cl [\hat{\partial} f]\big|_{M^{c}}$ holds. That is there exist sequences $x_i$ and $v_i$, with $v_i\in\hat{\partial} f(x_i)$ and $x_i\notin M$, such that $(x_i, f(x_i), v_i)$ converges to $(\bar{x},f(\bar{x}),\bar{v})$.
\end{lem}
%\begin{lem}[{\rm Accessibility}]\label{lem:main}
%Consider a lower semicontinuous, semi-algebraic function $f\colon\R^n\to\overline{\R}$ and a closed semi-algebraic set $M\subset\R^n$ such that $f\big|_M$ is finite. Assume $\dim [f]\big|_M < n$. Then any pair $(\bar{x},\bar{v})\in \R^n\times\R^n$, satisfying $\bar{x}\in M$ and $\bar{v}\in\partial f(\bar{x})$, can be approximated using points from outside of $M$. Namely, there exist sequences $x_i\in M^c$ and $v_i\in\hat{\partial} f(x_i)$ satisfying $(x_i, f(x_i), v_i)\to(\bar{x},f(\bar{x}),\bar{v})$.
%\end{lem} 
%Before presenting the proof, perhaps some comments are in order. Let's consider the special case of Lemma~\ref{lem:main} when $f(x)=\delta_Q(x)$, for some closed, convex, semi-algebraic set $Q$, and the closed semi-algebraic set $M$ is locally a manifold around some point $\bar{x}\in M$. For simplicity, assume that the dimension of the normal cone $N_Q(x)$ is constant for all $x\in M$. Then it is easy to see that condition $\dim [f]\big|_M < n$ is equivalent to requiring the strict inequality $\dim N_Q(\bar{x})<\dim N_M(\bar{x})$ to hold; that is, requiring that the set $Q$ not be ``sharp'' in some directions normal to $M$. We will discuss this notion in detail in Section[?]. Then the conclusion is that any normal vector $\bar{v}\in N_Q(\bar{x})$ may be accessed by vectors $v_i\in N_Q(x_i)$, for some points $x_i\notin M$ with $x_i\to\bar{x}$.     
{\pf
We first prove the lemma for the special case when $(\bar{x},f(\bar{x}),\bar{v})=(0,0,0)$. The general result will then easily follow. 
%We need the following simplifying result.
%\begin{claim}
%Without loss of generality, we can assume that the inclusion $\bar{v}\in\hat{\partial} f(\bar{x})$ holds. 
%\end{claim}
%{\pf Suppose that the lemma holds for such triples. We will show that this implies that the lemma holds in full generality. To this end, observe that
%by definition of the limiting subdifferential, there exists a sequence $(x_i,f(x_i),v_i)\in[\hat{\partial} f]$ converging to $(\bar{x},f(\bar{x}),\bar{v})$. Since $(\bar{x},f(\bar{x}),\bar{v})\notin\cl [\hat{\partial} f]\big|_{M^{c}}$, we deduce $x_i\in M$ for all large indices $i$. Hence $(x_i,f(x_i),v_i)\in[\hat{\partial} f]\big|_{M}$. Furthermore it is easy to check that for large indices $i$, we must have $(x_i,f(x_i),v_i)\notin\cl [\hat{\partial} f]\big|_{M^{c}}$. Hence the lemma applies to our sequence. Observe that for any real number $m>0$, we have 
%$$\inter\bigcup_{x\in M} \partial (f(\cdot) + \frac{1}{2}m|\cdot-x_i|^2)(x)=\inter\bigcup_{x\in M} \partial (f(\cdot) + \frac{1}{2}m|\cdot-\bar{x}|^2)(x) - m(x_i-\bar{x})$$. 

%}
Thus, assume that there exists a sequence of real number $m_i$ with $m_i\to\infty$, such that the inclusion 
\begin{equation}\label{eq:cond}
0\in \bd\bigcup_{x\in M} {m_i}x + \partial f(x), 
\end{equation}
holds. We must show that there exists a sequence $(x_i,f(x_i),v_i)\in[\hat{\partial} f]\big|_{M^{c}}$ converging to $(0,0,0)$.

We make some simplifying assumptions. Since $f$ is lower semicontinuous, there exists a real number $r>0$ such that $f\big|_{r\overline{\bf{B}}}\geq -1$. 
\begin{claim}
Without loss of generality, we can replace the function $f$ by $f_o:=f+\delta_{r\overline{\bf{B}}}$ and the set $M$ by $M_o:=M\cap \frac{1}{2}r\overline{\bf{B}}$. 
%assume $M\subset \frac{1}{2}r\overline{\bf{B}}$ and $f(x)=+\infty$ whenever $x\notin r\bf{\overline{B}}$. To see this, assume that the result holds under this assumption and consider $f_o:=f+\delta_{r\overline{\bf{B}}}$ and $M_o=M\cap \frac{1}{2}r\overline{\bf{B}}$. 
\end{claim}
{\pf Observe $(0,0,0)\in[\hat{\partial} f_o]\big|_{M_o}$.
Furthermore, we have 
$$[\partial f_o]\big|_{M_o}= [\partial (f+\delta_{r\overline{\bf{B}}})]\big|_{M\cap \frac{1}{2}r\overline{\bf{B}}} =  [\partial f]\big|_{M\cap \frac{1}{2}r\overline{\bf{B}}}\subset  [\partial f]\big|_{M}.$$ 
Combining this with (\ref{eq:cond}), we obtain $$0\in \bd\bigcup_{x\in M_o} m_ix +\partial f_o(x).$$
Consequently, if we replace the function $f$ by $f_0$ and the set $M$ by $M_0$, then the requirements of the lemma will still be satisfied. Now suppose that with this replacement, the result of the lemma holds. Then there exists a sequence $(x_i,f(x_i),v_i)\in[\hat{\partial}f_o]\big|_{M_o^{c}}$ converging to $(0,0,0)$. For indices $i$ satisfying $|x_i|< \frac{1}{2}r$, we have $x_i\notin M$ and $(x_i,f(x_i),v_i)\in[\hat{\partial}f]$. Thus restricting to large enough $i$, we obtain a sequence $(x_i,f(x_i),v_i)\in[\hat{\partial}f]\big|_{{M}^{c}}$ converging to $(0,0,0)$, as claimed. Therefore, without loss of generality, we can replace the function $f$ by $f_o$ and the set $M$ by $M_o$.}\qed \\

Thus to summarize, we have $$(\bar{x},f(\bar{x}),\bar{v})=(0,0,0),~~~f\big|_{r\overline{\bf{B}}}\geq -1,~~~ M\subset\frac{1}{2}r\overline{\bf{B}}, ~~~f(x)=+\infty ~{\rm for}~ x\notin r\overline{\bf{B}}.$$     

%We will now prove the lemma for the case when $\bar{v}$ is in the Frechet subdifferential $\hat{\partial} f(\bar{x})$, that is when the inclusion $(\bar{x},f(\bar{x}),\bar{v})\in[\hat{\partial}f]$ holds. The general case will follow by an easy limiting argument, which we state towards the end of the proof. 

%We will now prove the lemma for the case when $\bar{v}$ is in the Frechet subdifferential $\hat{\partial} f(\bar{x})$, that is when the inclusion $(0,0,0)\in[\hat{\partial}f]$ holds. The general case will follow by an easy limiting argument, which we state towards the end of the proof. 

We now define a certain auxiliary sequence of vectors $y_i$, which will allow us to construct the sequence $(x_i, f(x_i),v_i)$ that we seek. To this end, let $y_i$ be a sequence satisfying $y_i\to 0$ and 
\begin{equation}
y_i\notin \bigcup_{x\in M} m_i x+\partial f(x),  
\end{equation}
for each index $i$.
By (\ref{eq:cond}), such a sequence can easily be constructed. The motivation behind our choice of the sequence $y_i$ will soon become apparent.

The key idea now is to consider the following sequence of minimization problems.
%Since the set $D$ has dimension strictly less than $n$, it has an empty interior. So for any point $\bar{y}\in D$, we can choose a sequence $v_k\rightarrow \bar{y}$ with $v_k\notin D$. Now choose a sequence of real numbers $\lambda_k>0$ converging to zero such that $\lambda_k|v_k|$ converges to zero. Consider the sequence of minimization problems,
$$P(i):~~ \min_{x\in\R^n} ~\langle -y_i, x\rangle + m_i(d_M^2(x)+|x|^2)+f(x).$$ 
By compactness of the domain of $f$ and lower semi-continuity of $f$, we conclude that there exists a minimizer $x_i$ for the problem $P(i)$. For each index $i$, we have
\begin{align} 
\notag y_i\in\partial[m_i(d_M^2(\cdot)+|\cdot|^2)+f(\cdot)](x_i)&\subset\partial[m_i(d_M^2(\cdot)+|\cdot|^2)](x_i)+\partial f(x_i)\\
&\subset m_i(x_i-P_M(x_i)) +m_i x_i +\partial f(x_i),\label{eqn:inc}
\end{align}
where the inclusions follow from Proposition~\ref{prop:lip} and Theorem~\ref{thm:mar}.
We claim 
\begin{equation}
x_i\notin M,\label{eqn:notin}  
\end{equation}
for each index $i$. Indeed, if it were otherwise, from (\ref{eqn:inc}) we would have $$y_i\in m_i x_i +\partial f(x_i)\subset\bigcup_{x\in M} m_i x+\partial f(x),$$ thus contradicting our choice of the vector $y_i$. 

Now from (\ref{eqn:inc}), let $z_i\in P_M(x_i)$ be a vector satisfying 
\begin{equation}\label{eqn:defn}
v_i:=y_i -m_i(x_i-z_i)-m_i x_i\in \partial f(x_i).
\end{equation}
Our immediate goal is to show that the sequence $(x_i,f(x_i),v_i)\in [\partial f]\big|_{M^{c}}$ converges to $(0,0,0)$. To that end, evaluating the value function of $P(i)$ at $0$, we obtain 
$$0\geq \langle -y_i, x_i\rangle + m_i(d_M^2(x_i)+|x_i|^2)+f(x_i).$$ 
From (\ref{eqn:notin}), we deduce $x_i\neq 0$, and combining this with the inequality above, we obtain 
$$|y_i|\geq\langle y_i,\frac{x_i}{|x_i|}\rangle\geq m_i\frac{d_M^2(x_i)}{|x_i|}+m_i|x_i|+\frac{f(x_i)}{|x_i|}.$$

Since $y_i\to 0$, $m_i\to\infty$, and the function $f$ is bounded below, it is easy to see that $x_i$ converges to $0$. Furthermore, since we have $0\in\hat{\partial} f(0)$, we deduce $$\frac{f(x_i)}{|x_i|}\geq \frac{o(|x_i|)}{|x_i|}.$$ In particular, we conclude $m_i|x_i|\to 0$ and $f(x_i)\to 0$. Since $d_M(x_i)\leq |x_i|$, we deduce $m_i d_M(x_i)\to 0$. Hence from (\ref{eqn:defn}), we obtain $$|v_i|\leq |y_i| +m_i d_M(x_i)+m_i |x_i|\rightarrow 0.$$ Thus we have produced a sequence $(x_i,f(x_i),v_i)\in[\partial f]\big|_{M^{c}}$ converging to $(0,0,0)$ with $x_i\notin M$ for each index $i$. We are almost done. The trouble is that the vector $v_i$ is in the limiting subdifferential, rather than the Frechet subdifferential. However, this can be dealt with easily. Since $M$ is closed, it is easy to see that we can perturb the triples $(x_i,f(x_i),v_i)$, to obtain a sequence $(x'_i,f(x'_i),v'_i)\in[\hat{\partial}f]$ converging to $(0,0,0)$, still satisfying $x'_i\notin M$ for each index $i$. This completes the proof for the case when $(\bar{x},f(\bar{x}),\bar{v})= (0,0,0)$.\\

%Recall that, we have assumed that $\bar{v}$ is a Frechet subgradient, that is the inclusion $(0,0,0)\in[\hat{\partial}f]$ holds. We now  argue that the lemma holds when $0$ is a limiting subgradient. Thus assume $0\in\partial f(0)$. Hence there exists a sequence $(x_i,f(x_i),v_i)\in[\hat{\partial}f]$ converging to $(0,0,0)$. If there is a subsequence with $x_i\notin M$ for all large enough $i$, then we are done. Else for all large enough $i$, we have $x_i\in M$. Since we know that in this case, the lemma holds for each such triple $(x_i,f(x_i),v_i)$, the result easily follows.\\

Finally, we prove that the lemma holds when $(\bar{x},f(\bar{x}),\bar{v})\neq (0,0,0)$. Suppose that the point $(\bar{x},f(\bar{x}),\bar{v})$, the set $M$, and the function $f$ satisfy the requirements of the lemma. Now, consider the function $g(x):=f(x+\bar{x})-\langle\bar{v},x\rangle -f(\bar{x})$ and the set $N:=M-\bar{x}$. We will show that the function $g$, the set $N$, and the triple $(0,0,0)$ also satisfy the requirements of the lemma. To this end, observe $0\in N$ and $g(0)=0$. It is easy to verify the equivalence, $$v\in\hat{\partial} g(x)\Leftrightarrow v+\bar{v}\in\hat{\partial} f(x+\bar{x}).$$ Hence, clearly, $(0,0,0)\in[\hat{\partial} g]\big|_N$. Furthermore, the equation 
%To establish $(0,0,0)\notin\cl [\hat{\partial} g]\big|_{N^{c}}$, suppose otherwise. Then there exist sequences $x_i$ and $v_i$ with
%$$x_i\in N^{c}=M^{c}-\bar{x}, ~~v_i\in \hat{\partial} g(x_i), ~~(x_i,g(x_i),v_i)\rightarrow (0,0,0).$$ We deduce $x_i+\bar{x}\in M^{c}$ and $v_i+\bar{v}\in\hat{\partial} f(x_i+\bar{x})$. Since  $g(x_i)=f(x_i+\bar{x})-\langle\bar{v},x_i\rangle - f(\bar{x})$ converges to $0$, we obtain $f(x_i+\bar{x})\to f(\bar{x})$. Thus the triples $(x_i+\bar{x},f(x_i+\bar{x}),v_i+\bar{v})$ lie in $[\hat{\partial} f]\big|_{M^{c}}$ and converge to $(\bar{x},f(\bar{x}),\bar{v})$, which contradicts the assumption $(\bar{x},f(\bar{x}),\bar{v})\notin\cl[\partial f]\big|_{M^{c}}$. 
%Thus we can apply the lemma to the triple $(0,0,0)$, the function $g$, and the set $N$. We conclude that for all large enough real numbers $m$, the inclusion
$$\bigcup_{x\in N} \partial (g(\cdot) + \frac{1}{2}m|\cdot|^2)(x)= -\bar{v}+\bigcup_{x\in M} \partial (f(\cdot) + \frac{1}{2}m|\cdot-\bar{x}|^2)(x),$$
holds. Consequently, we deduce $0\in\bd \bigcup_{x\in N} \partial (g(\cdot) + \frac{1}{2}m|\cdot|^2)(x)$. We can now apply the lemma to the triple $(0,0,0)$, the function $g$, and the set $N$. Thus there exists a sequence $(x_i,f(x_i),v_i) \in[\hat{\partial} g]\big|_{N^{c}}$ with $(x_i, g(x_i), v_i)\to(0,0,0)$. Now observe that the sequence $(x_i+\bar{x},f(x_i+\bar{x}),v_i+\bar{v})$ lies in $[\hat{\partial} f]\big|_{M^c}$ and converges to $(\bar{x},f(\bar{x}),\bar{v})$, and hence the lemma follows.
}\qed 

In the semi-algebraic setting, Lemma~\ref{lem:main} yields the following important corollary. This corollary, in particular, will be crucial for proving our main result (Theorem~\ref{thm:local_dim}).
\begin{cor}\label{cor:main}
Consider a lower semicontinuous, semi-algebraic function $f\colon\R^n\to\overline{\R}$ and a closed semi-algebraic set $M\subset\R^n$ such that $f\big|_M$ is finite. Assume $\dim [\partial f]\big|_M < n$. Then any triple $(\bar{x},f(\bar{x}),\bar{v})$ in the restricted subjet $[\partial f]\big|_M$ can be accessed from the restricted subjet $[\hat{\partial} f]\big|_{M^{c}}$. That is, there exist sequences $x_i$ and $v_i$, with $v_i\in\hat{\partial} f(x_i)$ and $x_i\notin M$, such that $(x_i,f(x_i),v_i)\to(\bar{x},f(\bar{x}),\bar{v})$. Consequently, the inclusion $$[\partial f]\big|_M\subset \cl [\hat{\partial} f]\big|_{M^{c}},$$
holds. 
\end{cor} 
{\pf
Consider an arbitrary triple $(\bar{x},f(\bar{x}),\bar{v})\in[\hat{\partial} f]\big|_M$ and let $m$ be a positive real number.
Observe that the map $$\phi\colon [\partial f]\big|_M\to \gph(m(\cdot-\bar{x}) +\partial f(\cdot))\big|_M$$ $$(x,y,v)\mapsto (x,m(x-\bar{x})+v)$$ is bijective. Thus we deduce $$\dim\gph(m(\cdot-\bar{x}) +\partial f(\cdot))\big|_M=\dim [\partial f]\big|_M<n.$$ Hence, the set 
$$\bigcup_{x\in M} {m_i}(x-\bar{x}) + \partial f(x)= \bigcup_{x\in M} \partial (f(\cdot) + \frac{1}{2}m_i|\cdot-\bar{x}|^2)(x),$$ has dimension strictly less than $n$, and in particular has empty interior. Therefore, we have $$\bar{v}\in \bd\bigcup_{x\in M} \partial (f(\cdot) + \frac{1}{2}m_i|\cdot-\bar{x}|^2)(x).$$ By Lemma~\ref{lem:main}, we deduce that the inclusion $(\bar{x},f(\bar{x}),\bar{v})\in\cl [\hat{\partial} f]\big|_{M^{c}}$ holds. Consequently we obtain 
\begin{equation}\label{eq:te}
[\hat{\partial} f]\big|_M\subset\cl [\hat{\partial} f]\big|_{M^{c}}.
\end{equation}

Now consider a triple $(\bar{x},f(\bar{x}),\bar{v})\in[\partial f]\big|_M$. Then there exists a sequence $(x_i,f(x_i),v_i)\in [\hat{\partial} f]$ converging to $(\bar{x},f(\bar{x}),\bar{v})$. If there is a subsequence contained in $M^{c}$, then we are done. If not, then the whole sequence eventually lies in $M$, and then from (\ref{eq:te}) the result follows. 
}\qed
In order to prove our main result, we need to first establish a few simple propositions. We do so now.
\begin{prop}\label{prop:id}
Consider a semi-algebraic set $Q\subset\R^n$ and a point $\bar{x}\in Q$. Let $\{M_i\}$ be any stratification of $Q$. Then we have the identity $$\dim_Q(\bar{x})=\max_i\{\dim M_i: \bar{x}\in\cl M_i\}.$$
\end{prop}
{\pf Since there are finitely many strata, there exists some real number $\epsilon> 0$ such that for any $0<r<\epsilon$, we have  $$Q\cap B_{r}(\bar{x})=\bigcup_{i:\, \bar{x} \in {\mbox{{\scriptsize {\rm cl}}}}\, M_i} M_i\cap B_r(\bar{x}).$$
Hence, we deduce 
\begin{align*}
\dim (Q\cap B_{r}(\bar{x}))&= \max_i\{\dim (M_i\cap B_r(\bar{x})): \bar{x}\in\cl M_i\}\\
&= \max_i\{\dim M_i: \bar{x}\in\cl M_i\},
\end{align*}
where the last equality follows since the inclusion $\bar{x}\in\cl M_i$ implies that $M_i\cap B_r(\bar{x})$ is a nonempty open submanifold of $M_i$, and hence has the same dimension as $M_i$. Letting $r\to 0$ yields the result.
}\qed

\begin{defn}
{\rm Given a stratification $\{M_i\}$ of a semi-algebraic set $Q\subset\R^n$, we will say that a stratum $M$ is {\em maximal} if it is not contained in the closure of any other stratum.}
\end{defn}

\begin{rem}
{\rm Using the defining property of a stratification, we can equivalently say that given a stratification $\{M_i\}$ of a semi-algebraic set $Q\subset\R^n$, a stratum $M$ is maximal if and only if it is disjoint from the closure of any other stratum.
}
\end{rem}

\begin{prop}\label{prop:loc_max}
Consider a stratification $\{M_i\}$ of a semi-algebraic set $Q\subset\R^n$. Then given any point $\bar{x}\in Q$, there exists a maximal stratum $M$ satisfying $\bar{x}\in \cl M$ and $\dim M=\dim_Q(\bar{x})$. 
\end{prop}
{\pf By Proposition~\ref{prop:id}, we have the identity $$\dim_Q(\bar{x})=\max_i\{\dim M_i: \bar{x}\in\cl M_i\}.$$ Let $M$ be a stratum achieving this maximum. 
If there existed a stratum $M_i$ satisfying $M\subset\cl M_i$, then we would have $\dim M <\dim M_i$ and $\bar{x}\in \cl M\subset \cl M_i$, thus contradicting our choice of $M$. Therefore, we conclude that $M$ is maximal.
}\qed

We are now ready to prove the main result of this section.
\begin{thm}\label{thm:local_dim}
Let $f\colon\R^n\to\overline{\R}$ be a proper lower semicontinuous, semi-algebraic function. Then the subjet $[\hat{\partial}f]$ has local dimension $n$ around each of its points. The same holds for the limiting subjet $[\partial f]$.
\end{thm}
{\pf
We first prove the claim for the subjet $[\hat{\partial}f]$ and then the limiting case will easily follow.
Observe that the sets $\gph \hat{\partial} f$ and $[\hat{\partial}f]$ are in semi-algebraic bijective correspondence, via the map $(x,v)\longmapsto (x,f(x),v)$, and hence these two sets have the same dimension. Combining this observation with Theorem~\ref{thm:grd}, we deduce that the dimension of $[\hat{\partial}f]$ is $n$. Thus the local dimension of $[\hat{\partial}f]$ at any point is at most $n$. We must now establish the reverse inequality.  

Consider the subjet $[\hat{\partial}f]$ and the projection map $\pi\colon [\hat{\partial}f]\to \R^n$, which projects onto the first $n$ coordinates. Applying Theorem~\ref{thm:strat} to $\pi$, we obtain a finite partition of $[\hat{\partial}f]$ into disjoint semi-algebraic manifolds $\{{M_i}\}$ and a finite partition of the image $\pi([\hat{\partial}f])$ into disjoint semi-algebraic manifolds $\{L_j\}$, such that for each index $i$, we have $\pi(M_i)=L_j$ for some index $j$. 
%Observe that by the frontier condition, the local dimension mapping $\dim_{[f]}$ is constant along each $M_i$. 
%By Proposition~\ref{prop:id}, for each $(x,f(x),v)\in[f]$, we have $\dim_{[f]}(x,f(x),v)=\max\{\dim M_i|(x,f(x),v)\in \cl M_i\}.$  

Assume that the statement of the theorem does not hold. Thus there exists some point in the subjet $[\hat{\partial}f]$ at which $[\hat{\partial}f]$ has local dimension strictly less than $n$. Therefore, by Proposition~\ref{prop:loc_max}, there is a maximal stratum $M$ with $\dim M<n$. We now focus on this stratum.

\begin{lem}\label{lem:app}
$$\dim [\partial f]\big|_{\pi(M)}<n.$$
\end{lem}
{\pf
For each $x\in \pi(M)$, the set $M \cap \pi^{-1}(x)$ is open relative to $\pi^{-1}(x) $, since the alternative would contradict maximality of $M$. Thus $$\dim (M \cap\pi^{-1}(x))=\dim \pi^{-1}(x),$$ for each $x\in \pi(M)$. Therefore the sets $M$ and $[\hat{\partial}f]\big|_{\pi(M)}$, along with the projection map $\pi$, satisfy the assumptions of Proposition~\ref{prop:fiber}. Hence we deduce $\dim [\hat{\partial}f]\big|_{\pi(M)}= \dim M<n$. Observe $[\partial f]\setminus [\hat{\partial}f]\subset (\cl [\hat{\partial}f])\setminus [\hat{\partial}f]$. Hence as a direct consequence of Theorem~\ref{thm:grd}, we see $\dim ([\partial f]\setminus [\hat{\partial}f])\big|_{\pi(M)}\leq\dim ((\cl [\hat{\partial}f])\setminus [\hat{\partial}f]) <n$. Thus we conclude $\dim [\partial f]\big|_{\pi(M)} <n$, as was claimed.
}\qed

%Let $M$ be a stratum achieving $\max\{\dim M_i|(x,f(x),v)\in \cl M_i\}$. Observe $\dim M<n$ and, by maximality, the stratum $M$ is not contained in the closure of any other stratum. Hence in our notation, the stratum $M$ is locally maximal.
Let $U$ be a nonempty, relatively open subset of $\pi(M)$ such that $\cl U\subset \pi(M)$ and consider an arbitrary point $\bar{x}\in U$ with $(\bar{x},f(\bar{x}),\bar{v})\in M$.  Combining Corollary~\ref{cor:main} and Lemma~\ref{lem:app}, we conclude that there exists a sequence $(x_i,f(x_i),v_i)\in[\hat{\partial}f]$ converging to $(\bar{x},f(\bar{x}),\bar{v})$ where $x_i\notin \cl U$. Since $\bar{x}\in U$, we deduce $x_i\notin \pi(M)$ for all large enough $i$. Since there are finitely many strata, we conclude that the point $(\bar{x},f(\bar{x}),\bar{v})\in M$ is in the closure of some stratum other than $M$, thus contradicting maximality of $M$. Thus the subjet $[\hat{\partial}f]$ has local dimension $n$ around each of its points.      

Now for the limiting subjet, observe that for any real number $r>0$, we have $B_r(x,f(x),v)\cap [\hat{\partial}f]\neq\emptyset$. Hence it easily follows that $[\partial f]$ has local dimension $n$ around each of its points as well.
}\qed

%The following is a direct consequence of Theorem~\ref{thm:local_dim}.
%\begin{cor}\label{cor:lim}
%Let $f\colon\R^n\to\overline{\R}$ be a lower semicontinuous, semi-algebraic function. Then the subjet $[\partial f]$ has local dimension $n$ around each of its %points. 
%\end{cor}
%{\pf For any real number $r>0$, we have $B_r(x,f(x),v)\cap [\hat{\partial}f]\neq\emptyset$. Combining this fact with Theorem~\ref{thm:local_dim}, the result %follows. 
%}\qed

\begin{rem}\label{rem:fail}
{\rm If a semi-algebraic function $f\colon\R^n\to\overline{\R}$ is not lower semicontinuous, then the result of Theorem~\ref{thm:local_dim} can easily fail. For instance, consider the set $S:=\{x\in\R^2:|x|<1\}\cup\{(1,0)\}$. The local dimension of $[\partial \delta_S]$ at $((1,0),0,(1,0))$ is one, rather than two.  }
\end{rem}

\subsection{Geometry of the Clarke subdifferential mapping}\label{sec:Clarke}
Besides Frechet and limiting subdifferentials, there is another very important subdifferential, which we now define. In this subsection, we will restrict our attention to locally Lipschitz continuous functions. Recall that any locally Lipschitz continuous function $f\colon\R^n\to\R$ is differentiable almost everywhere, in the sense of Lebesgue measure.
\begin{defn}
{\rm Consider a locally Lipschitz function $f\colon\R^n\to\R$ and a point $x\in\R^n$. Let $\Omega\subset\R^n$ be the set of points where $f$ is differentiable. We define the {\em Clarke subdifferential} of $f$ at $x$ to be $$\partial_c f(x):=\conv\{\lim_{i\to\infty}\nabla f(x_i):x_i\to x, x_i\in\Omega\}.$$
}
\end{defn}
It is a nontrivial fact that for a locally Lipschitz continuous function $f\colon\R^n\to\R$ and a point $x\in\R^n$, we always have the equality $\partial_c f(x)=\conv\, \partial f(x)$. In particular, the inclusions $$\hat{\partial} f(x)\subset \partial f(x)\subset\partial_c f(x),$$  hold. Some interest in the Clarke subdifferential stems from the fact that this
subdifferential can be easier to approximate numerically. See for example \cite{BLO}. We should also note that the definition of the Clarke subdifferential can be extended to functions that are not locally Lipschitz continuous. Since we will not need this level of generality in this work, we do not pursue this further.

Consider a lower semicontinuous, semi-algebraic function $f\colon\R^n\to\R$. It is shown in \cite[Theorem 3.6]{dim} that the global dimension of the set $\gph \partial_c f$ is $n$. Since the Clarke subdifferential contains both the Frechet and the limiting subdifferentials, it is tempting to think that, just like in the Frechet and limiting cases, the graph of the Clarke subdifferential should have local dimension $n$ around each of its points. 

It can be shown that this indeed is the case when $n\leq 2$. In fact, this even holds for semi-linear functions for arbitrary $n$. (Semi-linear function are those functions whose domains can be decomposed into finitely many convex polyhedra so that the restriction of the function to each polyhedron is affine.) However for $n\geq 3$, as soon as we allow the function $f$ to have any curvature at all, the conjecture is decisively false. Consider the following illustrative example.

\begin{exa}\label{exa:Clarke}
{\rm
Consider the function $f\colon\R^3\to\R,$ defined by
\begin{displaymath}
   f(x,y,z) = \left\{
     \begin{array}{lr}
       \min\{x,y,z^2\} &, \mbox{\rm if}\,(x,y,z) \in \R_{+}^3\\
       \min\{-x,-y,z^2\} &,\mbox{\rm if}\, (x,y,z) \in \R_{-}^3\\
       0 &, \mbox{\rm{otherwise}.}       
     \end{array}
   \right.
\end{displaymath} 
It is standard to verify that $f$ is locally Lipschitz continuous and semi-algebraic. Let $\Gamma:=\conv\{(1,0,0),(0,1,0),(0,0,0)\}$. Consider the set of points $\Omega\subset\R^3$ where $f$ is differentiable. Then we have 
\begin{align*}
\conv\{\lim_{i\to\infty}\nabla f(\gamma_i):&\gamma_i\to (0,0,0), \gamma_i\in\Omega\cap\R_{+}^3\}=\\
&=\conv\{(1,0,0),(0,1,0),(0,0,0)\}=\Gamma, 
\end{align*}
and
\begin{align*} 
~~~~~~~\conv\{\lim_{i\to\infty}\nabla f(\gamma_i):&\gamma_i\to (0,0,0), \gamma_i\in\Omega\cap\R_{-}^3\}=\\
&=\conv\{(-1,0,0),(0,-1,0),(0,0,0)\}=-\Gamma.
\end{align*}
In particular, we deduce $\partial_c f(0,0,0)= \conv\{\Gamma\cup -\Gamma\}$. Hence the subdifferential $\partial_c f(0,0,0)$ has dimension two. 

Let $((x_i,y_i,z_i), v_i)\in\gph\partial_c f\big|_{\R_{+}^3}$ be a sequence converging to $((0,0,0),\bar{v})$, for some vector $\bar{v}\in\R^3$. Observe $v_i\in\conv\{(1,0,2z_i), (0,1,2z_i), (0,0,0)\}$. Hence, we must have $\bar{v}\in\Gamma$. Now consider a sequence $((x_i,y_i,z_i), v_i)\in\gph\partial_c f\big|_{\R_{-}^3}$ converging to $((0,0,0),\bar{v})$, for some vector $\bar{v}\in\R^3$. A similar argument as above yields the inclusion $\bar{v}\in -\Gamma$.  This implies that for any vector $\bar{v}$ in $\partial_c f(0,0,0)\setminus (\Gamma\cup -\Gamma)$, there does not exist a sequence $((x_i,y_i,z_i),v_i)\in\gph\partial_c f$ converging to $((0,0,0),\bar{v})$. Therefore for such a vector $\bar{v}$, there exists an open ball $B_\epsilon((0,0,0),\bar{v})$ such that $B_\epsilon((0,0,0),\bar{v})\cap \gph\partial_c f\subset \{(0,0,0)\}\times \partial_c f(0,0,0)$. Thus the local dimension of $\gph\partial_c f$ around the pair $((0,0,0), \bar{v})$ is two, instead of three.
}  
\end{exa}

\subsection{Composite optimization}
Consider a composite optimization problem $$\displaystyle\min_x\, g(F(x)),$$ where $g\colon\R^m\to\overline{\R}$ is a lower semicontinuous, semi-algebraic function and $F\colon\R^n\to\R^m$ is a smooth, semi-algebraic mapping. It is often computationally more convenient to replace the criticality condition $0\in\partial(g\circ F)(x)$ with the potentially different condition $0\in\nabla F(x)^{*}\partial g(F(x))$, related to the former condition by an appropriate chain rule. See for example the discussion of Lagrange multipliers \cite{lag}. Thus it is interesting to study the graph of the set-valued mapping $x\mapsto \nabla F(x)^{*}\partial g(F(x))$. In fact, it is shown in \cite[Theorem 5.3]{dim} that the dimension of the graph of this mapping is at most $n$. Furthermore, under some  assumptions, such as the set $F^{-1}(\dom \partial g)$ having a nonempty interior for example, this graph has dimension exactly $n$.     

In the spirit of our current work, we ask whether under reasonable conditions, the graph of the mapping $x\mapsto \nabla F(x)^{*}\partial g(F(x))$ has local dimension $n$ around each of its points. In fact, the answer is no. That is, subdifferential calculus does not preserve local dimension. As an illustration, consider the following example. 
\begin{exa}
{\rm Observe that for a lower semicontinuous function $f$, if we let $F(x)=(x,x)$ and $g(x,y)=f(x)+f(y)$, then we obtain $\nabla F(x)^{*}\partial g(F(x))=\partial f(x)+\partial f(x)$. 
Now let the function $f\colon\R\to\R$ be $f(x)=-|x|$. Then we have 
\begin{displaymath}
   \partial f(x) = \left\{
     \begin{array}{lr}
       1 & , x <0\\
       \{-1,1\} & , x =0\\
       -1& , x>0\\
     \end{array}
   \right.
\end{displaymath} 
The set $\gph \partial f$ has local dimension $1$ around each of its point, as is predicted by Theorem~\ref{thm:local_dim}. However, the graph of the mapping $x\mapsto\partial f(x)+\partial f(x)$ has an isolated point at $(0,0)$, and hence this graph has local dimension zero around this point, instead of one.

Furthermore, using Theorem~\ref{thm:local_dim}, we can now conclude that the mapping $x\mapsto\partial f(x)+\partial f(x)$ is not the subdifferential mapping of any semi-algebraic, lower semicontinuous function.  
}
\end{exa}   

\section{Consequences}
In this section, we present some consequences of Theorem~\ref{thm:local_dim}. Specifically, in Subsection~\ref{sub:minty} we develop a nonconvex, semi-algebraic analog of Minty's Theorem, and in Subsection~\ref{sub:sens} we derive certain sensitivity information about variational problems, using purely dimensional considerations. Both of these results illustrate that local dimension shows the promise of being a powerful, yet simple to use, tool in semi-algebraic optimization.  

\subsection{Analogue of Minty's Theorem}\label{sub:minty}
The celebrated theorem of Minty states that for a proper, lower semicontinuous, convex function $f\colon\R^n\to\overline{\R}$, the set $\gph \partial f$ is Lipschitz homeomorphic to $\R^n$ \cite{minty}. In fact, for each real number $\lambda>0$, the so called Minty map $(x,y)\mapsto \lambda x+y$ is such a homeomorphism. For nonconvex functions, Minty's theorem easily fails. However, one may ask if for a nonconvex, lower semicontinuous function $f\colon\R^n\to\overline{\R}$, a Minty type result holds locally around many of the points in the set $\gph \partial f$. In general, nothing like this can hold either. However, in the semi-algebraic setting, Theorem~\ref{thm:local_dim} does provide an affirmative answer. 

\begin{prop}\label{prop:exp}
If $Q\subset\R^p$ has local dimension $q$ around every point, then it is locally diffeomorphic to $\R^q$ around every point in a dense semi-algebraic subset.
%Let $f\colon\R^n\to\overline{\R}$ be a lower semicontinuous, semi-algebraic function. Then there exists a semi-algebraic set $D\subset \gph \partial f$, that is dense in $\gph \partial f$, and having the property that $\gph \partial f$ is locally diffeomorphic to $\R^n$ around any point $(x,v)\in D$. Analogous statement holds in the Frechet case. 
\end{prop}    
{\pf
Applying Theorem~\ref{thm:strat}, we obtain a stratification $\{M_i\}$ of $Q$. Let $D$ be the union of the maximal strata in the stratification. By Proposition~\ref{prop:loc_max}, we see that $D$ is dense in $Q$. Now consider an arbitrary point $x\in D$ and let $M$ be the maximal stratum containing this point. Since $Q$ has local dimension $q$ around $x$, we deduce that the manifold $M$ has dimension $q$. By maximality of $M$, there exists a real number $r>0$ such that $B_r(x)\cap Q=B_r(x)\cap M$, and hence $Q$ is locally diffeomorphic to $\R^n$ around $x$, as we claimed.
%Applying Theorem~\ref{thm:strat}, we obtain a stratification $\{M_i\}$ of $\gph\partial f$. Let $D$ be the union of the maximal strata in the stratification. By Proposition~\ref{prop:loc_max}, we see that $D$ is dense in $\gph\partial f$. Now consider an arbitrary point $(x,v)\in D$ and let $M$ be the maximal stratum containing this point. As a consequence of Corollary~\ref{cor:lim}, the manifold $M$ has dimension $n$. By maximality of $M$, there exists a real number $r>0$ such that $B_r(x,v)\cap \gph \partial f=B_r(x,v)\cap M$, and hence $\gph\partial f$ is locally diffeomorphic to $\R^n$ around $(x,v)$, as we claimed. The Frechet case is similar. 
}\qed

Consider a lower semicontinuous, semi-algebraic function $f\colon\R^n\to\overline{\R}$. Combining Proposition~\ref{prop:exp} and Theorem~\ref{thm:local_dim}, we see that $\gph\partial f$ is locally diffeomorphic to $\R^n$ around every point in a dense semi-algebraic subset. In fact, we can significantly strengthen Corollary~\ref{prop:exp}. Shortly, we will show that we can choose the local diffeomorphisms of Corollary~\ref{prop:exp} to have very simple form that is analogous to the Minty map.

We will say that a certain property holds for a {\em generic} vector $v\in\R^n$ if the set of vectors for which this property does not hold is a semi-algebraic set of dimension strictly less than $n$. In the semi-algebraic setting, this notion coincides with the measure-theoretic concept of ``almost everywhere''.  For a more in-depth discussion of generic properties in the semi-algebraic setting, see for example \cite{gen,dim}. 

\begin{defn}
{\rm
For a set $Q\subset\R^n$ and a map $\phi\colon Q\to\R^m$, we say that $\phi$ is {\em finite-to-one} if for every point $x\in\R^m$, the set $\phi^{-1}(x)$ consists of finitely many points.}
\end{defn}

We need the following proposition, which is essentially equivalent to \cite[Theorem 4.9,]{DM}. We sketch a proof below, for completeness.
\begin{prop}\label{prop1}
Let $Q\subset\R^n\times\R^n$ be a semi-algebraic set having dimension no greater than $n$. Then for a generic matrix $A\in\R^{n\times n}$, the map $$\phi_A\colon Q\to\R^n,$$ $$(x,y)\mapsto Ax+y,$$ is finite-to-one.  
\end{prop}
{\pf
Let $I\in\R^{n\times n}$ be the identity matrix and consider the matrix $[A,I]$. Let $L$ denote the nullspace of $[A,I]$. It is standard to check the equivalence 
\begin{equation}
Ax+y=b\Leftrightarrow\pi_{L^{\perp}}(x,y)=\pi_{L^{\perp}}(0,b)\label{eq:old}, 
\end{equation}
where $\pi_{L^{\perp}}$ denotes the orthogonal projection onto $L^{\perp}$. Recall that each element of a dense collection of $n$ dimensional subspaces of $\R^n\times\R^n$ can be written uniquely as $\rge [A, I]^{T}$ for some matrix $A$. From \cite[Theorem 4.9]{DM}, we have that for a generic $n$-dimensional subspace $U$ of $\R^n\times\R^n$, the orthogonal projection map $\pi_U\colon Q\to U$ is finite-to-one. Hence, we deduce that for a generic matrix $A\in\R^{n\times n}$, the corresponding projection map $\pi_{L^{\perp}}$ is finite-to-one. Combining this with (\ref{eq:old}), the result follows.
}\qed

\begin{prop}\label{prop2}
Consider a semi-algebraic set $Q\subset\R^n$ and a continuous, semi-algebraic function $p\colon Q\to\R^m$ that is finite-to-one. Then there exists a stratification of $Q$ such that for each stratum $M$, the map $p\big|_M$ is a diffeomorphism onto its image. 
\end{prop}
{\pf 
Applying Theorem~\ref{theorem:Hardt} to the map $p$, we obtain a partition of the image $p(Q)$ into semi-algebraic sets ${C_i}$ such that the map $p$ is semi-algebraically trivial over each $C_i$. Thus for each index $i$, and any point $c\in C_i$, there is a semi-algebraic homeomorphism $h\colon p^{-1}(C_i)\rightarrow C\times p^{-1}(c)$, such that the diagram, 
\begin{diagram}[height=1.7em]
p^{-1}(C_i) &\rTo^h   &C_i\times p^{-1}(c)\\
          &\rdTo_p  &\dTo_{\mbox{\scriptsize {\rm proj}}_{C_i}}\\
          &         &C_i
\end{diagram}
commutes.
 
Fix some index $i$. We will now show that the map $p$ is injective on any connected subset of $p^{-1}(C_i)$. To this effect, consider a connected subset $M\subset p^{-1}(C_i)$. Observe that the set $h(M)$ is connected. Since $p^{-1}(c)$ is a finite set, we deduce that there exists a point $v\in p^{-1}(c)$ such that the inclusion, 
\begin{equation}
h(M)\subset C_i\times \{v\}\label{eq: inc} 
\end{equation}
holds. Now given any two distinct points $x,y\in M$, since $h$ is a homeomorphism, we have $h(x)\neq h(y)$. Combining this with (\ref{eq: inc}), we deduce $p(x)={\rm proj}_{C_i}\circ h(x)\neq {\rm proj}_{C_i}\circ h(y)=p(y)$, as we needed to show.
   
Applying Theorem~\ref{thm:strat} to the map $p$, we obtain a finite partition of $Q$ into connected, semi-algebraic manifolds $\{M_i\}$ compatible with $\{p^{-1}(C_i)\}$, such that for each stratum $M_i$, the map $p\big|_{M_i}$ is smooth and $p$ has constant rank on $M_i$. Fix a stratum $M$. Since $M$ is connected, it follows from the argument above that $p$ is injective on $M$. Combining this observation with the fact that $p$ has constant rank on $M$, we deduce that $p\big|_{M}$ is a diffeomorphism onto its image. 
}\qed

We are now ready for the main result of this subsection.
\begin{thm}\label{thm:dif}
Consider a semi-algebraic set $Q\subset\R^{n\times n}$ that has local dimension $n$ around every point. Then for a generic matrix $A\in\R^{n\times n}$, the map $$\phi_A\colon Q\to\R^n,$$ $$(x,y)\mapsto Ax+y,$$ is a local diffeomorphism of $Q$ onto an open subset of $\R^n$, around every point in a dense semi-algebraic subset of $Q$.  
%Consider a semi-algebraic set Q\subset\R^{n\times n}. Then for a generic matrix $A\in\R^{n\times n}$, there exists a semi-algebraic set $D_A\subset Q$, that is dense in $Q$, and having the property that around any point $(\bar{x},\bar{y})\in Q$, the map $$\phi_A\colon Q\to\R^n,$$ $$(x,y)\mapsto Ax+y,$$ is a local diffeomorphism of $Q$ onto an open subset of $\R^n$. Analogous statement holds in the Frechet case. 
%Let $f\colon\R^n\to\overline{\R}$ be a lower semicontinuous, semi-algebraic function. Then for a generic matrix $A\in\R^{n\times n}$, there exists a semi-algebraic set $D_A\subset \gph \partial f$, that is dense in $\gph \partial f$, and having the property that around any point $(\bar{x},\bar{y})\in D_A$, the map $$\phi_A\colon \gph \partial f\to\R^n,$$ $$(x,y)\mapsto Ax+y,$$ is a local diffeomorphism of $\gph \partial f$ onto an open subset of $\R^n$. Analogous statement holds in the Frechet case. 
\end{thm}
{\pf By Proposition~\ref{prop1}, we have that for a generic matrix $A\in\R^{n\times n}$, the map $\phi_A$ is finite-to-one. Fix such a matrix $A$. Consider the stratification guaranteed to exist by applying Proposition~\ref{prop2} to the map $\phi_A$, and let $D_A$ be the union of the maximal strata in this stratification. By Proposition~\ref{prop:loc_max}, we see that $D_A$ is dense in $Q$. Consider a point $(\bar{x},\bar{y})\in D_A$, which is contained in some maximal stratum $M$. Since the set $Q$ has local dimension $n$ around each of its points, we deduce that the stratum $M$ is $n$-dimensional. Recall that the mapping $\phi_A\big|_M$ is a diffeomorphism onto its image. By maximality of $M$, there is a real number $\epsilon>0$ such that $B_{\epsilon}(\bar{x},\bar{y})\cap M=B_{\epsilon}(\bar{x},\bar{y})\cap Q$ and hence the restricted mapping $\phi_A\big|_{B_{\epsilon}(\bar{x},\bar{y})\cap Q}$ is a diffeomorphism onto its image.
 %$$\phi_A\big|_{B_{\epsilon}(\bar{x},\bar{y})\cap Q}\colon B_{\epsilon}(\bar{x},\bar{y})\cap Q\to \phi_A(B_{\epsilon}(\bar{x},\bar{y})\cap Q),$$ is a diffeomorphism. 
Consequently the image $\phi_A(B_{\epsilon}(\bar{x},\bar{y})\cap Q)$ is an $n$-dimensional submanifold of $\R^n$, and hence is an open subset of $\R^n$.  
}\qed

As a direct consequence of Theorem~\ref{thm:dif} and Theorem~\ref{thm:local_dim}, we obtain
\begin{cor}
Let $f\colon\R^n\to\overline{\R}$ be a lower semicontinuous, semi-algebraic function. Then for a generic matrix $A\in\R^{n\times n}$, the map $$\phi_A\colon \gph \partial f\to\R^n,$$ $$(x,y)\mapsto Ax+y,$$ is a local diffeomorphism of $\gph \partial f$ onto an open subset of $\R^n$ around every point in a dense semi-algebraic subset of $\gph \partial f$. Analogous statement holds in the Frechet case. 
\end{cor}

\subsection{Sensitivity}\label{sub:sens}

\begin{prop}\label{prop:pres}
Consider a semi-algebraic set $Q$ and a finite-to-one, continuous, semi-algebraic map $\phi\colon Q\to\R^m$. Then the map $\phi$ does not decrease local dimension, that is
$$\dim_Q(x)\leq\dim_{\rgel{\phi}}\phi(x),$$ for any point $x\in Q$. In particular, semi-algebraic homeomorphisms preserve local dimension. 
\end{prop}
{\pf 
By Proposition~\ref{prop2}, there exists a stratification of $Q$ into semi-algebraic manifolds $\{M_i\}$, such that for each maximal stratum $M$, the restriction $\phi\big|_M$ is a diffeomorphism onto its image. Fix some point $x\in Q$. By Proposition~\ref{prop:loc_max}, there is a maximal stratum $M$ satisfying $x \in \cl M$ and $\dim M=\dim_Q(x)$. Now since $\phi\big|_M$ is a diffeomorphism onto its image, we deduce that the manifold $\phi(M)$ has dimension $\dim_Q(x)$. By continuity of $\phi$, we have $\phi(x)\in\cl \phi(M)$. Hence, $$\dim_{\rgel\phi}\phi(x)\geq \dim \phi(M)=\dim_Q(x),$$ as we needed to show.   
%and consider the set $\rge\phi\cap B_{\epsilon}(\phi(x))$, where $\epsilon >0$ is some real number.
%We have $$\dim_{\rge\phi}\phi(x)\leq \dim \rge\phi\cap B_{\epsilon}(\phi(x))\leq \dim \phi^{-1}(B_{\epsilon}(\phi(x)))\leq \dim_Q(x),$$ where the last inequality holds because by continuity of $\phi$, the set $\phi^{-1}(B_{\epsilon}(\phi(x)))$ is a neighborhood of $\phi(x)$. Hence we must now establish the reverse inequality.
}\qed

\begin{prop}\label{prop:sens}
Let $Q\subset\R^n\times\R^n$ be a semi-algebraic set and suppose that $Q$ has local dimension $n$ at a point $(\bar{x},\bar{y})$. Consider the following parametric system, parametrized by matrices $A\in\R^{n\times n}$ and vectors $b\in\R^n$.
\begin{align*}
P(A,b):~~~~~ &(x,y)\in Q,\\
&Ax+y=b.
\end{align*}
Define the solution set, $S(A,b)$, to be the set of all pairs $(x,y)$ solving $P(A,b)$. Suppose that we have $(\bar{x},\bar{y})\in S(\bar{A},\bar{b})$, for some matrix $\bar{A}$ and vector $\bar{b}$. Fix some precision parameter $\epsilon >0$, and let $\Omega\subset\R^{n\times n}\times\R^n$ be the set of parameters $(A,b)$, for which the solution set $S(A,b)$ is finite and the intersection $S(A,b)\cap B_{\epsilon}(\bar{x},\bar{y})$ is nonempty.
Then for any real number $\delta>0$, the set $\Omega\cap B_{\delta}(\bar{A},\bar{b})$ has dimension $n^2+n$, and in particular has strictly positive measure.
\end{prop}
{\pf  
By Proposition~\ref{prop1}, for a generic matrix $A\in\R^{n\times n}$ the map $$\phi_A\colon Q\to\R^n,$$ $$(x,y)\mapsto Ax+y,$$ is finite-to-one. Denote this generic collection of matrices by $\Sigma$. Let $Q':=Q\cap B_{\epsilon}(\bar{x},\bar{y})$. Observe that for each matrix $A\in\Sigma$, the restriction $\phi_A\big|_{Q'}$ is still finite-to-one. For notational convenience, we will abuse notation slightly and we will always use the symbol $\phi_A$ to mean the restriction of $\phi_A$ to $Q'$, that is we now have $\phi_A\colon Q'\to\R^n$.

Fix some arbitrary real numbers $\delta,\gamma>0$, and let $N_{\delta,\gamma}(\bar{A},\bar{b}):=B_\delta(\bar{A})\times B_\gamma(\bar{b})$. We will show that the set $\Omega\cap N_{\delta,\gamma}(\bar{A},\bar{b})$ has dimension $n^2+n$. To this effect, observe that the inclusion, 
\begin{equation}\label{eq:weird}
\Omega\cap N_{\delta,\gamma}(\bar{A},\bar{b})\supset \{(A,b)\in \R^{n\times n}\times \R^n: A\in\Sigma\cap B_{\delta}(\bar{A}), b\in\rge\phi_A\cap B_{\gamma}(\bar{b})\},
\end{equation} 
holds. The set on the right hand side of (\ref{eq:weird}) is exactly the graph of the set-valued mapping,
$$F\colon\Sigma\cap B_{\delta}(\bar{A})\rightrightarrows \R^n,$$ $$A\mapsto \rge\phi_A\cap B_{\gamma}(\bar{b}).$$ Thus, in order to complete the proof, it is sufficient to show that $\gph F$ has dimension $n^2+n$. We will do this by showing that both the domain and the values of $F$ have large dimension.

First, we analyze the domain of $F$. Consider any matrix $A\in\Sigma\cap B_{\delta}(\bar{A})$. We have $$|\phi_A(\bar{x},\bar{y})-\bar{b}|= |(A\bar{x}+\bar{y})- (\bar{A}\bar{x}+\bar{y})|\leq |A-\bar{A}||\bar{x}|.$$
So by shrinking $\delta$, if necessary, we can assume $|\phi_A(\bar{x},\bar{y})-\bar{b}|<\gamma$.
Hence, we deduce
\begin{equation}\label{eqn:temp}
\phi_A(\bar{x},\bar{y})\in\rge\phi_A\cap B_{\gamma}(\bar{b}).
\end{equation}
In particular, we deduce that $F$ is nonempty valued on $\Sigma\cap B_{\delta}(\bar{A})$. Combining this with the fact that the set $\Sigma$ is generic, we obtain
\begin{equation}\label{eq:obv}
\dim \dom F=\dim \Sigma\cap B_{\delta}(\bar{A})=n^2. 
\end{equation}

We now analyze the set $F(A)$. Since the continuous map $\phi_A$ is finite-to-one and $Q'$ has local dimension $n$ at the point $(\bar{x},\bar{y})$, appealing to Proposition~\ref{prop:pres}, we obtain 
\begin{equation}\label{eqn:lc}
\dim_{\rgel\phi_A} {\phi_A(\bar{x},\bar{y})}=n.
\end{equation}
From (\ref{eqn:temp}) and (\ref{eqn:lc}), we obtain 
\begin{equation}\label{eq:last}
\dim F(A)=\dim \rge\phi_A\cap B_{\gamma}(\bar{b})=n,
\end{equation}
for all matrices $A\in\Sigma\cap B_{\delta}(\bar{A})$.
Finally combining (\ref{eq:obv}), (\ref{eq:last}), and Proposition~\ref{prop:const_gen}, we deduce $$\dim\gph F=\dim \dom F +\dim F(A)= n^2+n,$$ thus completing the proof.
}\qed

Thus we have the following corollary.
\begin{cor}\label{cor:sens}
Let $f\colon\R^n\to\overline{\R}$ be a lower semicontinuous, semi-algebraic function. Consider the following parametric system, parametrized by matrices $A\in\R^{n\times n}$ and vectors $b\in\R^n$.
\begin{align*}
P(A,b):~~~~~ &y\in \partial f(x),\\
&Ax+y=b.
\end{align*}
Define the solution set, $S(A,b)$, to be the set of all pairs $(x,y)$ solving $P(A,b)$. Suppose that we have $(\bar{x},\bar{y})\in S(\bar{A},\bar{b})$, for some matrix $\bar{A}$ and vector $\bar{b}$. Fix some precision parameter $\epsilon >0$, and let $\Omega\subset\R^{n\times n}\times\R^n$ be the set of parameters $(A,b)$, for which the solution set $S(A,b)$ is finite and the intersection $S(A,b)\cap B_{\epsilon}(\bar{x},\bar{y})$ is nonempty.
Then for any real number $\delta>0$, the set $\Omega\cap B_{\delta}(\bar{A},\bar{b})$ has dimension $n^2+n$, and in particular has strictly positive measure.
\end{cor}

To clarify Corollary~\ref{cor:sens}, consider a solution $(\bar{x},\bar{y})$ to the system $P(\bar{A},\bar{b})$. Then the content of Corollary~\ref{cor:sens} is that under small random (continuously distributed) perturbations to the pair $(\bar{A},\bar{b})$, with positive probability the perturbed system $P(A,b)$ has a strictly positive and finite number of solutions arbitrarily close to $(\bar{x},\bar{y})$.  \\

\noindent{\bf Acknowledgment}: Much of the current work has been done while the first and second authors were visiting CRM (Centra de Recerca Matem\`{a}tica) at Universitat Aut\`{o}nomo de Barcelona. The concerned authors would like to acknowledge the hosts for their hospitality. We thank Aris Daniilidis and J\'{e}r\^{o}me Bolte for fruitful discussions, and we also thank C.H. Jeffrey Pang for providing the illustrative Example~\ref{exa:Clarke}.

\bibliographystyle{plain}
\small
\parsep 0pt
\bibliography{dim_graph}

\end{document}